\documentclass{amsart}
\usepackage{amsopn,amstext,amsbsy,amsmath,amscd,amsthm,amsfonts}
\theoremstyle{plain}
\newtheorem{theorem}                 {\bf Theorem}     [section]
\newtheorem{proposition}  [theorem]  {\bf Proposition}
\newtheorem{corollary}    [theorem]  {\bf Corollary}
\newtheorem{lemma}        [theorem]  {\bf Lemma}
\newtheorem{satz}         [theorem]  {\bf Proposition}

\theoremstyle{definition}
\newtheorem{example}      [theorem]  {\bf Example}
\newtheorem{remark}       [theorem]  {\bf Remark}
\newtheorem{definition}   [theorem]  {\bf Definition}

\newcommand{\im}{\mathop{\rm im}}
\newcommand{\Hom}{\mathop{\rm Hom}\nolimits}
\newcommand{\Ind}{\mathop{\rm Ind}\nolimits}
\newcommand{\Coind}{\mathop{\rm Coind}\nolimits}

\numberwithin{equation}{section}

\allowdisplaybreaks

\author{Elias Kappos}

\title{$l^p$-cohomology for groups of type $FP_n$}



\address
{Mathematisches Institut, Georg-August-Universit\" at zu G\" ottingen, Bunsenstr. 3--5, D-37073 G\" ottingen, Germany}
\email{ekappos@uni-math.gwdg.de}

\begin{document}

\def\rn{\mathbb R}
\def\cn{\mathbb C}
\def\nn{\mathbb N}
\def\zn{\mathbb Z}
\def\zg{\zn[G]}
\def\zh{\zn[H]}
\def\lp{l^p(G)}
\def\lq{l^q(G)}
\def\lph{l^p(H)}
\def\hp{H^{(p)}_*(G)}
\def\hpgn{H^{(p)}_n(G)}
\def\hpgi{H^{(p)}_i(G)}
\def\hphi{H^{(p)}_i(H)}
\def\cpgn{H_{(p)}^n(G)}
\def\cpgi{H_{(p)}^i(G)}
\def\cphi{H_{(p)}^i(H)}
\def\cp{H_{(p)}^*(G)}
\def\norm#1{\parallel\! #1 \!\parallel}
\def\norme#1{\mid\! #1 \!\mid}

\begin{abstract}
Let $G$ be a group of type $FP_n$ and let $p>1$. In this paper we show that the reduced $l^p$-homology of $G$ is dual to the reduced $l^q$-cohomology for $\frac{1}{p}+\frac{1}{q}=1$. In our main theorem we show that for a group of type $FP_n$ with a central element of infinite order, the reduced $l^p$-cohomology vanishes. We generalize this fact for groups with infinitely many elements in the center of the group, for groups which are FCC, for groups with infinitely many finite conjugacy classes, for nilpotent groups, and for groups of polynomial growth.
\end{abstract}

\renewcommand{\thefootnote}{\fnsymbol{footnote}}

\maketitle

\section{Group Cohomological Preliminaries}

\footnotetext{The author was partially supported by DAAD and DFG} Everything in this section is standard in group cohomology, for example cf.\ \cite{brown}.
Let $G$ be a group and $\zg$ be the group ring over $\zn$, i.e. the free $\zn$-module generated by $G$. Then $\zn$ is a $\zg$-module with trivial $G$-action. A resolution of $\zn$ over $\zg$ is an exact sequence of $\zg$-modules
\begin{eqnarray*}
\dots\to F_2\xrightarrow{\partial_2}F_1\xrightarrow{\partial_1}F_0\xrightarrow{\varepsilon}\zn\to 0.
\end{eqnarray*}
Here the map $\varepsilon$ is the augmentation map. This resolution is called projective, if each $F_i$ is projective.
The standard resolution (which we will use only once in this paper), is of course, the most common, but the following theorem shows that it is not important which resolution we use for calculations. 
\renewcommand{\thefootnote}{\arabic{footnote}}

\begin{theorem}\label{theorem homotopieaequivalenz}
Let $F$ and $F'$ be projective resolutions of a $\zg$-module $M$, then there is an augmentation preserving chain map $f:F\to F'$, unique up to homotopy, and $f$ is a homotopy equivalence.
\end{theorem}

A resolution or partial resolution is said to be of {\it finite type}, if each $F_i$ is finitely generated. A group $G$ is said to be of type $FP_n$ if there exists a partial resolution s.t. $F_0,\dots,F_n$ are projective and finitely generated.

\begin{satz}\label{satz freie moduln}
Let $M$ be a $\zg$-module and $n\geq 0$. Then the following statements are equivalent:
\begin{itemize}
\item[i)] There exists a partial resolution $F_n\to\dots \to F_0\to M\to 0$, where the $F_i$ are free and finitely generated;
\item[ii)] $M$ is of type $FP_n$;
\item[iii)] $M$ is finitely generated and for every projective partial resolution 
\begin{equation*}
P_k\to\dots\to P_0\to M\to 0
\end{equation*}
of finite type with $k<n$, we have that $\ker\{P_k\to P_{k-1}\}$ is finitely generated.
\end{itemize}
\end{satz}

\begin{satz}\label{satz untergruppe typ fpn}
Let $H\subset G$ be a subgroup of finite index. Then $G$ is of type $FP_n$ if and only if $H$ is of type $FP_n$. 
\end{satz}

Let $G$ be a group, $F$ be a projective resolution of $\zn$ over $\zg$, and $M$ be a $\zg$-module. Then we define the {\it homology of $G$ with coefficients in $M$} by
\begin{eqnarray*}
H_*(G,M)=H_*(F\otimes_{\zg}M),
\end{eqnarray*}
and the {\it cohomology of $G$ with coefficients in $M$} by
\begin{eqnarray*}
H^*(G,M)=H^*(\Hom_{\zg}(F,M)).
\end{eqnarray*}

Finally, we introduce for a group $G$ and a subgroup $H$ the induction and co-induction for a $\zh$-module $M$ by
\begin{eqnarray*}
\Ind^G_HM=\zg\otimes_{\zh}M \text{\qquad and \qquad} \Coind^G_HM=\Hom_{\zh}(\zg,M).
\end{eqnarray*}

\begin{satz}{\bf (Shapiros Lemma)}\label{shapiro}
Let $H\subset G$ and let $M$ be a $\zh$-module. Then
\begin{eqnarray*}
H_*(H,M)\cong H_*(G,\Ind_H^G(M)) \text{\qquad and\qquad } H^*(H,M)\cong H^*(G,\Coind_H^G(M)).
\end{eqnarray*}
\end{satz}

\begin{lemma}\label{satz ind gleich coind}
Let $G$ be a group, $H$ be a subgroup of finite index $k$ and $M$ be a $\zh$-module. Then $\Ind_H^G M=\Coind_H^G M$. 
\end{lemma}

\section{Basic Definitions}

Let $\mathcal{F}(G)$ be the set of complex functions on the group $G$. Every $f\in\mathcal{F}$ can be represented as a formal sum $\sum\limits_{g\in G}f_gg$ with $f_g\in \cn$ and $f(g)=f_g$. 
For any real $p\geq 1$, $l^p(G)$ is the set of formal sums $f=\sum\limits_{g\in G}f_gg$ in $\mathcal{F}(G)$, where only countably many $f_g$ do not vanish and $\sum\limits_{g\in G}|f_g|^p=(|f|_p)^p<\infty$. Analogously we define $l^\infty(G)$ as those formal sums with only countable many $f_g\neq 0$ with $\sup_{g\in G}(f_g)=|f|_\infty<\infty$. Note that with this norm, $l^p(G)$ is a Banach space, and that $l^p(G)$ is a $\zg$-module because $G$ operates on $\lp$ (from both sides). 

\begin{definition}
Let $G$ be a group and $\dots\to F_2\xrightarrow{\partial_2} F_1\xrightarrow{\partial_1} F_0 \xrightarrow{\varepsilon}\zn\to 0$ be a projective resolution of $\zn$ over $\zg$. Then the $l^p$-(co-)homology of $G$ is defined as
\begin{itemize}
\item[] $\hpgn:=\ker(\partial_n\otimes 1)/\im(\partial_{n+1}\otimes 1)$,
\item[] $\cpgn:=\ker(\widetilde{\partial}_{n+1})/\im(\widetilde{\partial}_n)$,
\end{itemize}
where $\partial_i\otimes 1$ are the boundary maps in the complex $F_i\otimes_{\zg}\lp$ and $\widetilde{\partial_i}$ are the boundary maps in the complex $\Hom_{\zg}(F_i,\lp)$.
Finally, $\widetilde{\partial}_i$ is the dual map of $\partial_i$ defined by $\widetilde{\partial}_i(\varphi):=\varphi\circ\partial_i$.
In general, the image of the boundary operators are not closed in the kernel of the next boundary operator. Hence we have to define the {\it reduced $l^p$-(co-)homology of $G$} as
\begin{itemize}
\item[] $\overline{\hpgn}:=\ker(\partial_n\otimes 1)/\overline{\im(\partial_{n+1}\otimes 1)}$,
\item[] $\overline{\cpgn}:=\ker(\widetilde{\partial}_{n+1})/\overline{\im(\widetilde{\partial}_n)}$.
\end{itemize}
\end{definition}

A direct consequence of this definition is the fact that for a group $G$ and a $K(G,1)$-complex $X$ we have 
\begin{itemize}
\item[] $\hp=H_*(X;\lp)$ and
\item[] $\cp=H^*(X;\lp)$.
\end{itemize}

\begin{remark}
We need a norm on $\ker(\partial_n\otimes 1)$ and $\ker(\tilde{\partial}_{n+1})$ respectively for taking the closure in the definition of the reduced (co-)homology. But in this paper we are only working with groups of type $FP_n$ and for those the norm is natural.
\end{remark}

\section{Duality of homology and cohomology}

In the following we are interested in the topological duality of $l^p$-homology and $l^q$-cohomology, with the obvious requirement $\frac{1}{p}+\frac{1}{q}=1$. In the $l^2$-case this duality provides us with a standard tool. 
To prove our duality statement we need to review of some basic facts from functional analysis.

\begin{definition}\label{def megginson}
Let $X$ be a normed space and $A$ and $B$ subsets of $X$ and $X^*$, respectively. Then the {\it annihilators} $A^\perp$ and $^\perp B$ are defined by 
\begin{eqnarray*}
A^\perp:=\{x^*:x^*\in X^*, x^*x=0\;\forall x\in A\},
\end{eqnarray*} 
\begin{eqnarray*}
^\perp B:=\{x:x\in X, x^*x=0\;\forall x^*\in B\}.
\end{eqnarray*}
\end{definition}

\begin{proposition}\label{prop megginson}\cite[p. 93]{megginson}
Let $X$ be a normed space and let $A$ and $B$ be subsets of $X$ and $X^*$ respectively. Then $A^\perp$ and $^\perp B$ are closed subsets of $X^*$ and $X$ respectively.
Furthermore we have $^\perp(A^\perp)=\overline{A}$, if $A$ is a subset of $X$. 
\end{proposition}

\begin{theorem}\label{theorem megginson}\cite[pp. 94--95]{megginson}
Let $M$ be subset of a normed space $X$. Then there exists an isometric isomorphism that identifies $M^*$ with $X^*/M^\perp$. 

If $M$ is closed, then there exists an isometric isomorphism that identifies \linebreak $(X/M)^*$ with $M^\perp$. 
\end{theorem} 

\begin{lemma}\label{dualer kern}\cite[p. 289]{megginson}
Let $X$ and $Y$ be normed spaces and $T\in B(X,Y)$. Then $\ker(T)= {}^\perp(T^*(Y^*))$ and $\ker(T^*)=(T(X))^\perp$.
\end{lemma}

\begin{satz}\label{satz dual homolog}
Let $1<p\in\rn$, $n\in\nn$ and $G$ be a group of type $FP_n$. Then $\overline{\cpgi}$ is dual to $\overline{H^{(q)}_i(G)}$, with $\frac{1}{p}+\frac{1}{q}=1$, for $i\leq n$.
\end{satz}

\begin{proof}
The spaces $F_i\otimes_G\lp$ and $\Hom_G(F_i,\lq)$ are dual to each other, because the $F_i$ are projective and finitely generated for $i\leq n$: 
\begin{eqnarray*}
(F_i\otimes_G\lp)^*&\cong&F_i^*\otimes_G\lp^*\\&\cong&Hom_G(F_i^{**},l^p(G)^*)\\&\cong&Hom_G(F_i,l^q(G)).
\end{eqnarray*}
In a similar way the boundary operators in the chain complexes are dual to each other, because $-\otimes_G\lp$ is a covariant functor and $Hom_G(-,\lq)$ is a contravariant functor (cf. \cite[Chapter III.1]{brown}). 

For simplicity we write $C_i:=F_i\otimes_G\lp$, with boundary operator \linebreak $\varphi_i:C_{i+1}\to C_i$, and analogously $C_i^*:=Hom_G(F_i,\lq)$ with boundary operator $\varphi^*_i:C^*_i\to C^*_{i+1}$.  
The sequence
\begin{eqnarray*}
0\to\{f\in C^*_i:f|_{\ker(\varphi_{i-1})}\equiv 0\}\to C^*_i\to \ker(\varphi_{i-1})^*\to 0
\end{eqnarray*}
is exact because $C^*_i$ is the direct sum of the other two spaces. Restricting this sequence to $\{f\in C^*_i:f|_{\overline{\im(\varphi_i)}}\equiv 0\}$ proves the exactness of the following sequence:
\begin{eqnarray*}
0\to\{f\in C^*_i:f|_{\ker(\varphi_{i-1})}\equiv 0\}&\to& \{f\in C^*_i:f|_{\overline{\im(\varphi_i)}}\equiv 0\} \\*  &\to& \{f\in(\ker(\varphi_{i-1})^*:f|_{\overline{\im(\varphi_i)}}\equiv 0\}\to 0.
\end{eqnarray*} 
The basic facts of functional analysis provide us with isometric identifications
\begin{eqnarray*}
\ker(\varphi_i^*)
&\overset{\ref{dualer kern}}{=}&\im(\varphi_{i})^\perp\\
&\overset{\ref{def megginson}}{=}&\{f\in C^*_i:f|_{\im(\varphi_i)}\equiv 0\}\\
&=&\{f\in C^*_i:f|_{\overline{\im(\varphi_i)}}\equiv 0\}.
\end{eqnarray*}
The last equation holds because of the continuity of the functions $f$.
We get similarly
\begin{eqnarray*}
\overline{\im(\varphi_{i-1}^*)}
&\overset{\ref{prop megginson}}{=}&(^\perp\im(\varphi_{i-1}^*))^\perp\\
&\overset{\ref{dualer kern}}{=}&\ker(\varphi_{i-1})^\perp\\
&\overset{\ref{def megginson}}{=}&\{f\in C^*_i:f|_{\ker(\varphi_{i-1})}\equiv 0\}.
\end{eqnarray*}
Furthermore we have
$\left(\left.\ker(\varphi_{i-1})\right/\overline{\im(\varphi_i)}\right)^*=\{f\in(\ker(\varphi_{i-1}))^*:f|_{\overline{\im(\varphi_i)}}\equiv 0\}$, by Theorem \ref{theorem megginson}.
\end{proof}

\begin{remark}\label{auswertungspaarung}
The proof implies that a pairing on the complexes induces a pairing on the (co-)homology groups. In particular this is true for the evaluation pairing. 
\end{remark}

\section{First Examples}

\begin{example}\label{triviales beispiel}
Let $p\geq 1$ and let $G$ be the trivial group. Then $\hpgi=0$ for $i\geq 1$ and $H_0^{(p)}(G)=\cn$. 
\end{example}

\begin{proof}
The trivial group has the resolution $0\to\zn\xrightarrow{id}\zn\to 0$. The $l^p$-functions on $G$ are isomorphic to $\cn$, hence the tensor product of the resolution with the $l^p$-functions gives us $0\to\cn$.
\end{proof}

\begin{satz}\label{satz endlicher index}
Let $G$ be a group and  $H\subset G$ be a subgroup of finite index of $G$. Then $\hpgi=\hphi$ and $\cpgi=\cphi$, for all $i\in\nn_0$ and $p\geq 1$.
\end{satz}

\begin{proof}
This is a direct consequence of Proposition \ref{shapiro}, Lemma \ref{satz ind gleich coind} and the identifications
\begin{eqnarray*}
Ind_H^G(\lph)=\zg\otimes_{\zn [H]}\lph\cong\lp.
\end{eqnarray*}
\end{proof}

\begin{corollary}\label{beispiel endliche gruppe}
Let $p\geq 1$ and let $G$ be a finite group. Then $\hpgi=0$ for $i\geq 1$ and $H_0^{(p)}(G)=\cn$. 
\end{corollary}

\begin{proof}
Every finite group has the trivial group as a subgroup of finite index. Thus the statement follows from Example \ref{triviales beispiel} and Proposition \ref{satz endlicher index}.
\end{proof}

\begin{lemma}\label{lemma zyklische gruppe}
Let $G$ be an infinite cyclic group with generator $t$. Then $G$ has vanishing reduced $l^p$-homology and vanishing reduced $l^p$-cohomology, for all $1<p<\infty$.
\end{lemma}

\begin{proof}
In this case we have a resolution
\begin{eqnarray*}
0\to\zg\xrightarrow{t-1}\zg\to\zn\to 0.
\end{eqnarray*}
Hence $\hp$ is the homology of
\begin{eqnarray*}
\dots\to 0\to \lp\xrightarrow{t-1}\lp
\end{eqnarray*}
and $\cp$ is the cohomology of
\begin{eqnarray*}
\lp\xrightarrow{t-1}\lp\to 0\to\dots.
\end{eqnarray*}
Thus $\overline{H_1^{(p)}(G)}=\overline{H^0_{(p)}(G)}=0$ because the Kernel of the map $t-1$ is $0$. Theorem \ref{satz dual homolog} proves the rest of the statement.
\end{proof} 

\begin{remark}
Lemma \ref{lemma zyklische gruppe} is a good illustration of the difference between reduced and unreduced homology and cohomology. The image of the map $t-1$ is not equal to $\lp$. The element $t$ itself for example is not in the image of this map. Only the closure of the image repairs this problem. 
\end{remark}

\section{A Vanishing Theorem}

In \cite{puls} Puls asserts a vanishing result for first $l^p$-cohomology. However there appears to be a serious gap in the proof. In \cite{martin} Martin and Valette provided a different proof by generalizing his result for a bigger class of groups. In this section we provide another proof, using other methods, which yields a much stronger result. 

\begin{lemma}\label{basislemma}
Let $G$ be a group and $\underline{h}\in \mbox{Z}(G)$ be an element in the center of $G$. Let $\varepsilon :F\to \zn$ be a projective resolution of $\zn$ over $\zg$. Then $h$, multiplication by the element $\underline{h}$, is a homotopy equivalence, and chain homotopic to the identity on $F^*:=Hom_G(F,\lp)_*$.
\end{lemma}

\begin{proof}
This proof has three parts. First, we have to show that $h$ is a chain map and a chain equivalence; second, we have to prove the lemma for the standard resolution; and finally, we have to prove the lemma for any resolution.

Part 1: $h$ is bounded because $|\varphi|_p=|\underline{h}\varphi|_p=|h(\varphi)|_p$. $\underline{h}$ is an element in the center of $G$, therefore the equality $\varphi(x)=(a\varphi)(a^{-1}x)$ remains valid under multiplication by $\underline{h}$. 
Multiplication by an element of $G$ is a $\zg$-module homomorphism and the boundary operators are $\zg$-homomorphisms, too; hence they commute with multiplication by an element of $\zg$. It follows that $h$ is a chain map. To get a chain equivalence we need the map $h^{-1}$, but this is multiplication by $\underline{h}^{-1}$.

Part 2: Let $\varepsilon:F\to\zn$ be the standard resolution of $\zn$ over $\zg$, i.e. the boundary operator $\partial$ is defined as 
\begin{eqnarray*}
\partial\varphi(x_0,\dots ,x_n)=\sum_{i=0}^n(-1)^i\varphi(x_0,\dots ,x_{i-1},\widehat{x_i},x_{i+1},\dots ,x_n)
\end{eqnarray*}
where  $\;\widehat{}\;$  has the usual meaning (to omit the distinguished term). Define \linebreak $j:F_n\to F_{n+1}$ as
\begin{eqnarray*}
j\varphi(x_0,\dots ,x_n)=\sum_{k=0}^n(-1)^{k+1}\varphi(x_0,\dots ,x_k,h(x_k),\dots ,h(x_n)),
\end{eqnarray*}
then $j$ is again a $\zg$-module homomorphism and bounded. Hence 
\begin{eqnarray*}
(\partial j + j\partial)(\varphi(x))
&=&\sum_{k=0}^n\sum_{i=0}^{k}(-1)^{i+k+1}\varphi(x_0,\dots ,\widehat{x_i},\dots ,x_k,h(x_k),\dots ,h(x_n))\\
&&+\sum_{k=0}^n\sum_{i=k}^{n}(-1)^{i+k+2}\varphi(x_0,\dots ,x_k,h(x_k),\dots ,\widehat{h(x_i)},\dots ,h(x_n))\\
&&+\sum_{i=0}^n\sum_{k=0}^{i-1}(-1)^{i+k+1}\varphi(x_0,\dots ,x_k,h(x_k),\dots ,\widehat{h(x_i)},\dots ,h(x_n))\\
&&+\sum_{i=0}^n\sum_{k=i+1}^{n}(-1)^{i+k}\varphi(x_0,\dots ,\widehat{x_i},\dots ,x_k,h(x_k),\dots ,h(x_n))\\
&=&\sum_{k=0}^n\sum_{i=0}^{k}(-1)^{i+k+1}\varphi(x_0,\dots ,\widehat{x_i},\dots ,x_k,h(x_k),\dots ,h(x_n))\\
&&+\sum_{k=0}^n\sum_{i=0}^{k-1}(-1)^{i+k}\varphi(x_0,\dots ,\widehat{x_i},\dots ,x_k,h(x_k),\dots ,h(x_n))\\
&&+\sum_{k=0}^n\sum_{i=k}^{n}(-1)^{i+k+2}\varphi(x_0,\dots ,x_k,h(x_k),\dots ,\widehat{h(x_i)},\dots ,h(x_n))\\
&&+\sum_{k=0}^n\sum_{i=k+1}^{n}(-1)^{i+k+1}\varphi(x_0,\dots ,x_k,h(x_k),\dots ,\widehat{h(x_i)},\dots ,h(x_n))\\
&=&\sum_{k=0}^n(-1)^{2k+1}\varphi(x_0,\dots,\widehat{x_k},h(x_k),dots,h(x_n))\\
&&+\sum_{k=0}^n(-1)^{2k+2}\varphi(x_0,\dots,x_k,\widehat{h(x_k)},dots,h(x_n))\\
&=&\varphi(x_0,\dots ,x_n)-\varphi(h(x_0),\dots ,h(x_n))\\
&=&\varphi(x)-h(\varphi(x)),
\end{eqnarray*}
with $x=(x_0,\dots,x_n)$. Therefore $j$ is a chain homotopy between the identity and $h$. 

Part 3:
Let $\varepsilon' :F'\to\zn$ be another projective resolution of $\zn$ over $\zg$. By Theorem \ref{theorem homotopieaequivalenz} there exists a homotopy equivalence $f:C^*\to F^*$. 
\end{proof}

\begin{remark}\label{bemerkung Gruppenring uebertragung}
The proof of Lemma \ref{basislemma} can be transferred almost verbatim to the case $\underline{h}\in\zg^\times\cap \mbox{Z}(\zg)$. The only necessary modification is to show that the new map $h$ is bounded. This follows by the H\"older inequality.
\end{remark}

\begin{remark}
A similar lemma can be formulated for $F_*=F\otimes_G\lp$. The proof is also very similar.
\end{remark}

Martin and Valette \cite{martin} have done the following calculations for the first $l^p$-cohomology. As mentioned earlier, our methods extend their result to higher (co-)homology groups. 

\begin{theorem}\label{satz unendlich viele zentrale gruppenelemente}
Let $1<p\in\rn$, $n\in\nn$, and $G$ be a group of type $FP_n$ with infinitely many elements in the center $Z(G)$ of $G$. Then $\overline{\hpgi}=0$ and $\overline{\cpgi}=0$ for $i\leq n$. 
\end{theorem}

\begin{proof}
By Proposition \ref{satz freie moduln}, there exists a partial resolution of $\zn$ over $\zg$ of free, finitely generated $\zg$-modules $F_i$. This means $F_i\otimes_G\lp$ is isomorphic to $\lp^{m_i}$ and $Hom_G(F_i,\lq)$ to $\lq^{m_i}$, respectively. Here $m_i$ is the dimension of $F_i$. 

Now we assume that there exists a $0\neq\overline{x}\in\overline{\hpgi}$ represented by 
\begin{eqnarray*}
x=\left(\sum_{h\in G}\xi_{1,h}h,\dots,\sum_{h\in G}\xi_{m_i,h}h\right).
\end{eqnarray*}

Let $\overline{b}$ be the induced pairing discussed in Remark \ref{auswertungspaarung}. This is well-definied because the images of the boundary operators are $\zg$-submodules.  

By duality, cf.\ Proposition \ref{satz dual homolog}, there exists a $0\neq\overline{y}\in\overline{H_{(q)}^i(G)}$ represented by 
\begin{eqnarray*}
y=\left(\sum_{h\in G}\zeta_{1,h}h,\dots,\sum_{h\in G}\zeta_{m_i,h}h\right)
\end{eqnarray*}
such that $\overline{b}(\overline{y},\overline{x})\neq 0$. In particular there are only countably many $\zeta_i$ and $\xi_i$ not equal to $0$. The induced evaluation pairing of $\overline{x}$ and $\overline{y}$ is given by
\begin{eqnarray*}
\overline{b}(\overline{y},\overline{x})=\sum_{j=1}^{m_i}\sum_{h\in G}\zeta_{j,h}\xi_{j,h}.
\end{eqnarray*}
This sum converges absolutely as it is a finite sum of absolutely convergent series.

Now let $D$ be an infinite countable subset of $Z(G)$ such that $1\in D$; if $g_i\in D$ then $g_i^{-1}\in D$, and if $g_i,g_j\in D$ then $g_ig_j\in D$. Let the elements in $D$ be indexed by the integers in such a way that $1=g_0$. Then
\begin{eqnarray*}
\overline{b}(\overline{y},g_i\overline{x})=\sum_{j=1}^{m_i}\sum_{h\in G}\zeta_{j,h}\xi_{j,g_i^{-1}h}
\end{eqnarray*}
converges for all $i\in\zn$. 

Let $\varepsilon>0$. Then there exist finite sets $A$ and $B$ in $G$ and a subset $I$ of $\{1,\dots,m_i\}$ such that $\norm{x}_p^p-\sum\limits_{a\in I\times A} |\xi_a|^p<\varepsilon$ and $\norm{y}_q^q-\sum\limits_{b\in I\times B} |\zeta_{b}|^p<\varepsilon$. Hence there exist countable sets $A'=\{cg_i\mid c\in A, i\in\zn\}$ and $B'=\{dg_i\mid d\in B, i\in\zn\}$ containing the elements of the generalized orbits $x^{D}$ of elements of $A$ and $B$, respectively (confer \cite{huppert}).  With the short hands $\tilde{A}=G^{m_i}\backslash (I\times A')$ and $\tilde{B}=G^{m_i}\backslash (I\times B')$ we have:
\begin{eqnarray*}
b(y,x)&=&\sum_{j=1}^{m_i}\sum_{h\in G}\zeta_{j,h}\xi_{j,h}\\
&=&\sum_{\substack{c\in \tilde{A},\\ d\in \tilde{B}}}\zeta_d\xi_c
+\sum_{\substack{c\in \tilde{A},\\ d\in (I\times B')}}\zeta_d\xi_c
+\sum_{\substack{c\in (I\times A'),\\d\in \tilde{B}}}\zeta_d\xi_c
+\sum_{\substack{c\in (I\times A'),\\ d\in (I\times B')}}\zeta_d\xi_c\\
\end{eqnarray*}
where $i=i'$ for the pairs $c=(i,c')$ and $d=(i',d')$ in the same summand. 
Taking norms and using the H\"older inequality, the first three summands are bounded by $\varepsilon(\varepsilon+\norm{x}_p+\norm{y}_q)$. Translation by $g_i$ does not change these summands. 

So we are only interested in the last one. This summand can be rearrenged to look like
\begin{eqnarray}\label{formel in satz unendlich}
\sum_{i\in I}\sum_{\substack{a\in \bar{A},\\b\in \bar{B}}}\sum_{r=-\infty}^\infty \zeta_{i,g_rb}\xi_{i,g_ra},
\end{eqnarray}
where $\bar{A}$ (resp. $\bar{B}$) is a subset of $A$ (resp. $B$) which contains only one representative of every orbit which is in $A'$ (resp. $B'$).

Translation by $g_i$ gives us 
\begin{eqnarray*}
\sum_{i\in I}\sum_{\substack{a\in \bar{A},\\b\in \bar{B}}}\sum_{r=-\infty}^\infty \zeta_{i,g_rb}\xi_{i,g_i^{-1}g_ra}.
\end{eqnarray*}

We can now use a similiar cut off argument on each orbit, because $D$ has infinitely many elements. To be precise, we get a finite subset $R\subset\zn$ with 

\begin{eqnarray*}
\sum\limits_{i\in I}\sum\limits_{r=-\infty}^\infty|\xi_{i,g_ra}|^p-\sum\limits_{i\in I}\sum\limits_{r\in R}|\xi_{i,g_ra}|^p<\varepsilon
\end{eqnarray*}
for all $a\in\bar{A}$ and a finite subset $S\subset\zn$ with 
\begin{eqnarray*}
\sum\limits_{i\in I}\sum\limits_{r=-\infty}^\infty|\zeta_{i,g_rb}|^q-\sum\limits_{i\in I}\sum\limits_{r\in S}|\zeta_{i,g_rb}|^q<\varepsilon
\end{eqnarray*} 
for all $b\in\bar{B}$. Set $T=R\cup S$.

Splitting up of the sum again, we get 
\begin{eqnarray*}
\sum_{i\in I}\sum_{\substack{a\in \bar{A},\\b\in \bar{B}}}\sum_{r=-\infty}^\infty \zeta_{i,g_rb}\xi_{i,g_ra}=\sum\limits_{i\in I}\sum_{\substack{a\in \bar{A},\\b\in \bar{B}}}\sum_{r\notin T}\zeta_{i,bg_r}\xi_{i,ag_r}+\sum\limits_{i\in I}\sum_{\substack{a\in \bar{A},\\b\in \bar{B}}}\sum_{r\in T}\zeta_{i,bg_r}\xi_{i,ag_r}
\end{eqnarray*}

Taking the norm as before, the first term is bounded by $\varepsilon\norm{x}_p$, if we translate by $g_i$. 

Choosing now $i$ such that for $g_j=g_i^{-1}g_r$ the index $j$ is not in $T$ for all $r\in T$, 
we get for all $\varepsilon>0$ a $i_\varepsilon \in\zn$ with $b(y,g_{i_\varepsilon}x)<W\cdot\varepsilon$ and $W$ constant. 

Hence in contradiction to our assumption there is for every $x$ representing $\overline{x}\in\overline{\hpgi}$ and every $y$ representing $\overline{y}\in\overline{H_{(q)}^i(G)}$ a net $\{g_{i_\varepsilon}\}$ with $\lim\overline{b}(\overline{y},g_{i_\varepsilon}\overline{x})=0$, and by Lemma \ref{basislemma}, we have that the multiplication with every $g_i$ is a homotopy equivalence. 
\end{proof}

The following corollary generalizes work of Puls, as mentioned in the beginning of this section. 

\begin{corollary}\label{satz zentrales Element unendlicher Ordnung}
Let $1<p\in\rn$, $n\in\nn$, and $G$ be a group of type $FP_n$ with a central element of infinite order. Then $\overline{\hpgi}=\overline{\cpgi}=0$ for $i\leq n$.
\end{corollary}

\section{Generalizations of the Vanishing Theorem}

\begin{corollary}
Let $1<p\in\rn$, $n\in\nn$, and $G$ be an infinite group of type $FP_n$, and let this group be {\it FCC}, i.e. all conjugacy classes have only finitely many elements. Then $\overline{\hpgi}=0$ and $\overline{\cpgi}=0$ for $i\leq n$. 
\end{corollary}

\begin{proof}
Let $g$  be an element in $G$ and $<\!\!g\!\!>$ be the conjugacy class of $g$ in $G$. Furthermore let $\mbox{Z}(g)$ be the centralizer  of $g$ and $\mbox{Z}(G)$ be the center of $G$. Then
\begin{eqnarray*}
\norme{\,<\!\!g\!\!>\,}=\norme{G/\mbox{Z}(g)}.
\end{eqnarray*}
$G$ is of type $FP_n$, so $G$ is in particular finitely generated. Let now $\{g_1,\dots,g_k\}$ be a finite system of generators of $G$. Hence $\mbox{Z}(G)=\bigcap_{i=1}^k\mbox{Z}(g_i)$. Therefore
\begin{eqnarray*}
\norme{G/\mbox{Z}(G)}\leq \prod_{i=1}^k\norme{G/\mbox{Z}(g_i)}<\infty.
\end{eqnarray*}
that is, the center of $G$ has finite index. Thus there are infinite many elements in the center of the group. Hence the statement follows from Proposition \ref{satz unendlich viele zentrale gruppenelemente}.
\end{proof}

Before we come to the next generalization, we have to take a look at the structure of the group ring $\zg$. Let $k$ be a finite conjugacy class and $\hat{k}$ the sum of the elements of $k$ as an element of $\zg$.

\begin{lemma}\cite{passman1}\label{lemma fcc}
The center $\mbox{Z}(\zg)$ of the group ring $\zg$ has the set \linebreak $\{\hat{k}:k \text{ is a finite conjugacy class in }G\}$ of the sums of all finite conjugacy classes as a $\zn$-base.
\end{lemma}

\begin{lemma}\label{basislemma2}
Let $G$ be a group and let $\underline{h}\in \mbox{Z}(\zg)\setminus\{0\}$ be an element in the center of $\zg$ not equal to $0$. Let $\varepsilon :F\to \zn$ be a projective resolution of $\zn$ over $\zg$. Then $h$, multiplication by the element $\underline{h}$, is a homotopy equivalence and chain homotopic to the identity on $F^*:=Hom_G(F,\lp)_*$.
\end{lemma}

\begin{proof}
Let $k$ be a conjugacy class in $G$ with $n_k$ elements. Then $\underline{b_k}=\sum_{g\in k}g$ is an element in the basis of the center of the group ring by Lemma \ref{lemma fcc}. Let $\varphi$ be an element in $F^i$. $\underline{b_k}$ is an element in the center of the group ring, thus we have 
\begin{eqnarray*}
\sum_{g\in k} g\varphi(x)=
(\underline{b_k}\varphi)(x)=\varphi(\underline{b_k}x)=\sum_{g\in k}\varphi (gx).
\end{eqnarray*}
So, multiplication by $\underline{b_k}$ from the left equals multiplication from the right, up to a permutation of the terms in the sum. This means that essentially we multiply by something similar to the sum of $n_k$ elements of the center of the group. So we deduce from Lemma \ref{basislemma} that our multiplication is chain homotopic to $n_k\cdot\mbox{id}$, and hence chain homotopic to the identity.

The $b_k$ form a basis of the center of the group ring, thus every element $\underline{h}\in \mbox{Z}(\zg)\setminus\{0\}$ can be written as $\sum\lambda_k\underline{b_k}$ with at least one $\lambda_k$ not equal to zero. Hence multiplication by $\underline{h}$ is chain homotopic to $\sum\norme{\lambda_k}\cdot\mbox{id}=n_{h}\mbox{id}$. But this is again chain homotopic to the identity.
\end{proof}

\begin{theorem}\label{theorem fcc}
Let $1<p\in\rn$, $n\in\nn$, and $G$ be a group of type $FP_n$ with infinitely many finite conjugacy classes. Then $\overline{\hpgi}=0$ and $\overline{\cpgi}=0$ for $i\leq n$. 
\end{theorem}

\begin{proof}
By Lemma \ref{lemma fcc} the center of the group ring has infinitely many different elements. Hence we have to show that we can replace an element of the center of the group by an element of the center of the group ring in the proof of Theorem \ref{satz unendlich viele zentrale gruppenelemente}.

The generalized definition of an orbit (cf.\ \cite{huppert}) can be generalized a second time by considering the orbits of a set of ring elements. Let
\begin{eqnarray*}
S=\{h_r=\sum\limits_{j=1}^{n_r}(\alpha_r)_{g_j}g_j|h_r\in \mbox{Z}(\zg)\} 
\end{eqnarray*}
be a countable subset of the basis of the center of the group ring given by Lemma \ref{lemma fcc}. 
In particular this means that the $\alpha_r$ can only be $1$ or $0$, and furthermore the group elements $g_j$ can have a coefficient not equal to zero only in one $h_r$. Without loss of generality we can write every $h_r$ as $h_r=\sum\limits_{j=1}^{n_r}g_{r_j}$.  To get an orbit similar to the one in Proposition \ref{satz unendlich viele zentrale gruppenelemente} we need a multiplicative inverse of every central element in the group ring. By the proof of Lemma \ref{basislemma2}, we can repair the missing inverse by a function $g_r$ such that $g_rh_r\simeq h_rg_r\simeq \mbox{id}$. Taking this into account and reviewing the proof of Theorem \ref{satz zentrales Element unendlicher Ordnung}, we see that using these new orbits does not change anything up to equation \ref{formel in satz unendlich}. The only difference is an excess of indices. 
The term in equation \ref{formel in satz unendlich} has to be reformulated because we do not multiply anymore by $g_r$ but by $h_r$, an element of the group ring. In other words, we get another summation:
\begin{eqnarray*}
\sum_{i\in I}\sum_{\substack{a\in \bar{A},\\b\in \bar{B}}}\sum_{r=1}^\infty\sum_{j=1}^{n_r} \zeta_{i,g_{r_j}b}\xi_{i,g_{r_j}a}.
\end{eqnarray*}
So in this case we need two more cut-off arguments instead of one. The first one of these chooses the maximum of the relevant $j$ for all $r$, such that the inner sum only goes up top this maximum and the limit of the summation becomes independent of $r$. 
To be precise there are only finitely many indices $u$ and $v$ such that $\norm{y}^p_p-\sum\zeta_{u}<\varepsilon$ and $\norm{x}^p_p-\sum\xi_{v}<\varepsilon$. We define $n$ to be the maximum of the occuring $j$, and if $n>n_r$, we take  for all $n_r<j<n$ that $g_{r_j}$ is some $g$ such that $\zeta_{(i,g)}$ and $\xi_{(i,g)}$ small enough. 

Hence we can change the inner two summation signs. The second cut-off argument is analogous to the last cut-off argument in the proof of Theorem \ref{satz unendlich viele zentrale gruppenelemente}.
Together with Lemma \ref{basislemma2} this proves the statement.
\end{proof}

\begin{remark}
In this proof, we have used orbits of the form $\{x\cdot h_i\mid i\in I, h_i\in S \}$. It is possible to prove the statement by ``orbits" of the form $\{x\cdot h_1\cdot\dots\cdot h_r\mid h_i\in S\}$. This has the benefit that we do not need inverses and it simplifies the indices a little bit. However, the proof given is very similar to the proof of Theorem \ref{satz zentrales Element unendlicher Ordnung}, which is why we took that route.
\end{remark}

\begin{satz}\label{satz nilpotent}
Let $1<p\in\rn$ and let $G$ be a finitely generated, infinite nilpotent group (cf.\ \cite{robinson}). Then $\overline{\hpgi}=0$ and $\overline{\cpgi}=0$ for all $i\in\nn$. 
\end{satz}

\begin{proof}
By \cite[p. 213]{brown} every finitely generated, nilpotent group $G$ is of type $FP_\infty$.
If now $G$ has an infinite center, the statement follows from Proposition \ref{satz unendlich viele zentrale gruppenelemente}. If not, there exists a central series and a $k\in\nn$ such that $G_k$ has infinitely many elements and $G_{k-1}$ has only finitely many. By the definition of a nilpotent group $G_k/G_{k-1}$ is in the center of $G/G_{k-1}$. Thus we have, with $g\in G$ and $g_k\in G_k$ that $g^{-1}g_kg=g^{-1}gg_k=g_k$ modulo $G_{k-1}$. Hence the conjugacy class of $g_k$ is contained in the coset $g_kG_{k-1}$. The conjugacy class is finite because $G_{k-1}$ has only finitely many elements. Hence we have infinitely many elements with finite conjugacy class in $G$, because $g_k$ was arbitrary in $G_k$. The statement follows now from Theorem \ref{theorem fcc}. 
\end{proof}

\begin{corollary}
Let $1<p\in\rn$ and let $G$ be a finitely generated, infinite group of polynomial growth. Then $\overline{\hpgi}=0$ and $\overline{\cpgi}=0$ for all $i\in\nn$. 
\end{corollary}

\begin{proof}
Gromov \cite{gromov2} tells us that every finitely generated group with polynomial growth has a nilpotent subgroup of finite index. Proposition \ref{satz nilpotent} together with Proposition \ref{satz endlicher index} imply the statement. 
\end{proof}

\end{document}